\documentstyle[fleqn,12pt]{article}
\begin{document}

\def\ch{{\rm{ch}}\,}
\def\Ker{{\rm{Ker}}\,}


 \title{Bosonization of \\
        Admissible Representations of
        $U_q(\widehat{sl}_2)$ \\
        at Level $-\frac{1}{2}$ and $q$-Vertex Operators}

 \author{Yoshitaka Koyama \\ \\
 Research Institute for Mathematical Sciences, \\
 Kyoto University, Kyoto 606, Japan}

 \date{}
 \maketitle

 \begin{abstract}

  We construct a representation
  of $U_q(\widehat{sl}_2)$ at level $-1/2$
  by using the bosonic Fock spaces.
  The irreducible modules are obtained
  as the kernel of a certain operator,
  in contrast to the construction by Feingold and Frenkel
  for $q=1$ where such a procedure is not necessary.
  We also bosonize the $q$-vertex operators
  associated with the vector representation.

 \end{abstract}


 \section{Introduction}

  The $q$-vertex operators introduced by Frenkel and Reshetikhin(\cite{fr})
  have important roles in the theory of solvable lattice models.
  In particular, at a positive integral level,
  the $q$-vertex operators are closely related to
  generalizations of the XXZ spin-chain
  and the RSOS models.(\cite{dfjmn},\cite{djo},\cite{iijmnt})

  In \cite{kw}, Kac and Wakimoto introduced
  a class of irreducible highest weight representations
  of affine Lie algebras called admissible representations
  which is wider than the class of integrable representations.
  Admissible representations have a fractional level,
  and therefore they are not integrable in general.
  They proved the character formula for admissible representations,
  which is a generalization of the Weyl-Kac character formula.
  From this formula,
  we see that
  the set of the characters of
  admissible representations of the same level
  is closed under the modular transformations.
  Furthermore, the minimal series of the Virasoro algebra
  are constructed from the admissible representations
  by tensoring the level 1 integrable representation(the coset construction).
  This shows that the concept of the admissible representations
  is a natural generalization of that of the integrable representations.

  Feingold and Frenkel constructed
  level $-1/2$ highest weight admissible representations
  of $\widehat{sp}_{2n}$ in terms of bosons(\cite{ff}).
  In $\widehat{sp}_2=\widehat{sl}_2$ case,
  the representations are constructed as follows.
  Let $\{ \beta_n, \gamma_n | n \in I \}$ ($I={\bf Z}$ or ${\bf Z}+1/2$)
  be the following bosonic $\beta$-$\gamma$ system 
  $$ [\beta_m , \beta_n]=[\gamma_m , \gamma_n]=0, \quad
     [\beta_m , \gamma_n]=\delta_{m+n,0} , $$
  and let ${\cal F}$ be the Fock space generated by the vacuum vector
   $|vac \rangle$ satisfying
  $$ \beta_m |vac \rangle = \gamma_n |vac \rangle =0 \quad
     {\rm for} \; m \geq 0, n > 0 . $$
  We set
  $$ E(n)=-\frac{1}{2}\sum_{n=m+m'}  \beta_m  \beta_{m'}   \; , $$
  $$ H(n)=-           \sum_{n=m+m'} :\beta_m  \gamma_{m'}: \; , $$
  $$ F(n)= \frac{1}{2}\sum_{n=m+m'}  \gamma_m \gamma_{m'}  \; , $$
  where $: \quad :$ is the normal ordering defined by
  $$ : \beta_m \gamma_n :
          = \left \{ \begin{array}{ll}
                       \beta_m \gamma_n  & (m<n) \\
                       \frac{1}{2}(\beta_m \gamma_n +\gamma_n \beta_m)
                       & (m=n) \\
                       \gamma_n \beta_m  & (m>n)
                       \end{array} \right.                       $$
  Then, $\{ E(n),H(n),F(n)|n \in {\bf Z} \}$
  forms $\widehat{sl}_2$ of level $-1/2$
  and ${\cal F}$ is completely reducible
  as an $\widehat{sl}_2$-module.
  The decomposition of ${\cal F}$
  into the irreducible components
  is calculated by characters(\cite{lu}).
  In \cite{ff}, the decomposition was proved
  by using the action of the Virasoro algebra
  via the Sugawara construction.

  Let $\phi_1(z)$, $\phi_2(z)$ be two independent bosonic fields
  normalized by
  $$ \phi_1(z) \phi_1(w) \sim -\log(z-w) , $$
  $$ \phi_2(z) \phi_2(w) \sim  \log(z-w) . $$
  It is known that 
  the bosonic $\beta$-$\gamma$ system can be
  expressed by $\phi_1(z)$ and $\phi_2(z)$ as follows(\cite{fms}).
  $$ \beta(z)=\sum_{n \in I} \beta_n z^{-n-1}
     \longrightarrow
     :\partial_z \phi_2(z) e^{\phi_1(z)+\phi_2(z)}: \; , $$
  $$ \gamma(z)=\sum_{n \in I} \gamma_n z^{-n}
     \longrightarrow
     :e^{-\phi_1(z)-\phi_2(z)}: \; . $$
  The Fock space of the bosonic $\beta$-$\gamma$ system
  is obtained as $\Ker Q^-$, where $Q^-=\oint :e^{-\phi_2(z)}:$. 
  Using the above picture, we can rewrite the representation
  by Feingold and Frenkel in terms of $\phi_1(z)$, $\phi_2(z)$.

  The aim of this paper is
  to construct a $q$-analogue of the representation
  in terms of $\phi_1(z)$, $\phi_2(z)$
  and bosonize the corresponding $q$-vertex operators
  associated with the vector representation (cf.\cite{m}).
  In Section 2,
  we review the known results about the quantum affine algebras
  and admissible representations.
  In Section 3,
  we construct the $q$-analogue of the above representation.
  In Section 4, we bosonize the $q$-vertex operators.
  Section 5 is devoted to the proofs
  which are not given in the previous sections.


 \section{Preliminaries}

   In this section, we recall
   the definition of the quantum affine algebra $U_q(\widehat{sl}_2)$
   in order to fix the notations
   and review the results about
   the admissible representations in \cite{kw}.

   Let $P$ and $P^*$ be free ${\bf Z}$-modules as follows.
   $$ P:=
      {\bf Z} \Lambda_0 \oplus {\bf Z} \Lambda_1 \oplus {\bf Z} \delta
      \; , \quad
      P^*:=
         {\rm Hom}(P,{\bf Z})
        ={\bf Z} h_0 \oplus {\bf Z} h_1 \oplus {\bf Z} d \; . $$
   We call $P$ and $P^*$
   the weight lattice and the coroot lattice respectively.
   The pairing is given by
   $ \langle \Lambda_i , h_j \rangle =\delta_{ij} $,
   $ \langle \Lambda_i , d \rangle =0             $,
   $ \langle \delta , h_j \rangle =0              $,
   $ \langle \delta , d \rangle =1                $.
   Let
   $ \alpha_0
        = 2 \Lambda_0 - 2 \Lambda_1 + \delta $,
   $ \alpha_1
        = 2 \Lambda_1 - 2 \Lambda_0 $
   be the simple roots.
   We fix an invariant symmetric bilinear form on $P$ such that
   $(\alpha_i|\alpha_i)=2$, $(\alpha_0|\alpha_1)=-2$.

   We denote the base field ${\bf Q}(q^{\frac{1}{4}})$ by ${\bf K}$.
   The quantum affine algebra $U_q(\widehat{sl}_2)$
   is the ${\bf K}$-algebra generated by the symbols
   $ \{ \, \gamma^{\pm \frac{1}{2}}, \, K, \, q^{\pm d},
        \,  h_k, \, x^{\pm}_l \,
        ( k \in {\bf Z} \setminus \{ 0 \}, \, l \in {\bf Z} ) \} $
   which satisfy the following defining relations(\cite{d}).
   \[ \gamma^{\pm\frac{1}{2}}
      \in {\it Center \; of} \; U_q(\widehat{sl}_2) \; ,
      \quad
      \gamma^{\frac{1}{2}} \gamma^{-\frac{1}{2}}=1 \; ,
      \quad
      K K^{-1} = K^{-1} K = 1 \; , \]
   $$ q^d q^{-d} = q^{-d} q^d = 1 \; ,
      \quad
      q^d x^{\pm}_k q^{-d} = q^k x^{\pm}_k \; ,
      \quad
      q^d h_k q^{-d} = q^k h_k \; ,
      \quad
      q^d K = K q^k \; , $$
   $$ [ h_k , K ] = 0 \; ,
      \quad 
      [ h_k , h_l ]
              =\delta_{k+l,0} \frac{[2k]}{k}
                              \frac{\gamma^k - \gamma^{-k}}{q-q^{-1}} \; ,
       \quad
       [ h_k , x^{\pm}_l ]
                =\pm \frac{[2k]}{k}
                     \gamma^{\mp \frac{|k|}{2}} x^{\pm}_{k+l} \; , $$
   $$ K x^{\pm}_k K^{-1}=q^{\pm 2} x^{\pm}_k \; ,
      \quad
      x^{\pm}_{k+1} x^{\pm}_{l} - q^{\pm 2} x^{\pm}_{l} x^{\pm}_{k+1}
       =q^{\pm 2} x^{\pm}_k x^{\pm}_{l+1} - x^{\pm}_{l+1} x^{\pm}_k \; , $$
   $$ [ x^+_k , x^-_l ] = \frac{1}{q-q^{-1}}
                 (\gamma^{(k-l)/2} \psi_{k+l}
                             -\gamma^{(l-k)/2} \varphi_{k+l}) $$
   where
     \[ \sum^{\infty}_{k=0} \psi_k z^{-k}
                     =K \, {\rm exp} \left( (q-q^{-1})
                             \sum^{\infty}_{k=1} a_k z^{-k}
                                                       \right) \; , \]
     \[ \sum^{\infty}_{k=0} \varphi_{-k} z^{-k}
                     =K^{-1} {\rm exp} \left( -(q-q^{-1})
                                  \sum^{\infty}_{k=1} a_{-k} z^{-k}
                                                       \right) \; . \]
    Here we used the standard notation
    $$ [n]=\frac{q^n-q^{-n}}{q-q^{-1}} . $$
    The Chevalley generators $\{ e_i, f_i, t_i \}$ are given by
    $$ t_0 = \gamma K^{-1} , \;
       e_0 = x^-_1 K^{-1}  , \;
       f_0 = K_1 x^+_{-1}  , \;
       t_1 = K ,     \;
       e_1 = x^+_0 , \;
       f_1 = x^-_0 . $$
    The algebra $U_q(\widehat{sl}_2)$ has a Hopf algebra structure
    with the coproduct $\Delta$ and the antipode $S$ given as follows.
    $$ \Delta (t_i) = t_i \otimes t_i
       \; , \;
       \Delta (q^d) = q^d \otimes q^d
       \; , $$
    $$ \Delta (e_i) = e_i \otimes 1 + t_i \otimes e_i
       \; , \;
       \Delta (f_i) = f_i \otimes t_i^{-1} + 1 \otimes f_i \; , $$
    $$ S(t_i) = t_i^{-1}
       \; , \;
       S(q^d) = q^{-d}
       \; , \;
       S(e_i) = -t_i^{-1}e_i
       \; , \;
       S(f_i) = -f_i t_i \; . $$
   Throughout this paper, we denote
   $U_q(\widehat{sl}_2)$ by $U_q$
   and the subalgebra of $U_q$ generated by
   $\{ t_i^{\pm 1},e_i,f_i \, (i=0,1) \}$ by $U'_q$.
   We denote by $V(\lambda)$
   the irreducible highest weight
   $U_q$(or $U'_q$)-module with highest weight $\lambda$.
   We fix a highest weight vector of $V(\lambda)$
   and denote it by $|\lambda \rangle$.
   If $W_i$ $(i=1,2)$ has a weight decomposition
   $W_i=\oplus_{\mu} W_{i,\mu}$,
   we define their completed tensor product by
   $$ W_1 \widehat{\otimes} W_2
      =\bigoplus_{\mu}
      (\prod_{\mu=\mu_1+\mu_2} W_{1,\mu_1} \otimes W_{2,\mu_2}) . $$

   Admissible representations are
   irreducible highest weight representations
   whose highest weight is an admissible weight(\cite{kw}).
   Admissible weights are defined as follows.
   \\
   {\bf Definition.(\cite{kw},\cite{kw2})}
   {\it
   We call a weight $\lambda$ admissible
   if it satisfies the following two conditions.
   \[ \hspace{-.2in}
      1) \; \;
         \langle \lambda + \rho , \alpha^{\vee} \, \rangle
         \notin \{ 0, -1, -2, -3, \cdots \}
         \quad
         for \; all \; real \; positive \; coroots \; \alpha^{\vee} , \]
   \[ \hspace{-.2in}
      2) \; \;
         {\bf Q} R^{\lambda} = {\bf Q} \Pi^{\vee} \]
   where
   $\rho$ is the sum of all fundamental weights,
   $\Pi^{\vee}$ is the set of simple coroots and }
   $$    R^{\lambda}
         = \{ \alpha^{\vee} : a \; positive \; real \; coroot
             | \langle \lambda + \rho , \alpha^{\vee} \, \rangle
               \in {\bf Z} \} . $$
   \hspace{\fill} $\Box$
   \\
   The classification of admissible weights is given in \cite{kw2}.
   In the $A^{(1)}_1$ case,
   the set of admissible weights of level
   $m=t/u$ ($t,u \in {\bf Z}$, $(t,u)=1$, $u \geq 1$)
   is given by
   \[ \, \{ (m-n+k(m+2)) \Lambda_0 + (n-k(m+2)) \Lambda_1 \]
   $$ | \; 0 \leq n \leq t+2u-2, \; 0 \leq k \leq u-1 \;
        (n,k \in {\bf Z}) \, \} . $$
   In this paper,
   we concentrate on the level $-1/2$ case.
   The level $-1/2$ admissible weights are as follows.
   $$ \lambda_1 = - \frac{1}{2} \Lambda_0 , \;
      \lambda_2 = - \frac{3}{2} \Lambda_0 + \Lambda_1 , \;
      \lambda_3 = - \frac{1}{2} \Lambda_1 , \;
      \lambda_4 = \Lambda_0 - \frac{3}{2} \Lambda_1 . $$
   In the level $-1/2$ case,
   the characters of the irreducible highest weight $\widehat{sl}_2$-modules
   $L(\lambda)$ have the following formulas due to \cite{lu}.
   $$ \ch L(\lambda_1) + p^{\frac{1}{2}} \ch L(\lambda_2)
      = \frac{e^{-\frac{1}{2} \Lambda_0}}
             {(p^{\frac{1}{2}} z^{ \frac{1}{2}};p)_{\infty}
              (p^{\frac{1}{2}} z^{-\frac{1}{2}};p)_{\infty}} \; , $$
   $$ \ch L(\lambda_3) + \ch L(\lambda_4)
      = \frac{e^{-\frac{1}{2} \Lambda_1}}
             {(p z^{-\frac{1}{2}};p)_{\infty}
              (  z^{ \frac{1}{2}};p)_{\infty}} \; , $$
   where
   $p=e^{-\delta}$, $z=e^{-\alpha_1}$
   and
   $(a;p)_{\infty}=\prod^{\infty}_{n=0}(1-a p^n)$.

   We use the difference operator
   $_{q^{\frac{1}{2}}} \partial_z$ defined by
   $$ _{q^{\frac{1}{2}}} \partial_z f(z)
     :=\frac{f(q^{\frac{1}{2}}z)-f(q^{-\frac{1}{2}}z)}
            {(q^{\frac{1}{2}}-q^{-\frac{1}{2}})z} . $$


 \section{Construction of representation}

   Let $\{ a_n, b_n | n \in {\bf Z} \}$ be a set of operators
   satisfying the following commutation relations for $n \neq 0$.
   $$ [ a_n, a_{-n} ] = \frac{[2n][-\frac{1}{2}n]}{n}
      \quad , \quad
      [ b_n, b_{-n} ] = n . $$
   The other commutation relations are zero.
   We define the Fock module ${\cal F}_{l_1,l_2}$
   by the defining relations
   $$ a_n | l_1, l_2 \rangle = 0 \quad (n>0)
      \; , \quad
      b_n | l_1, l_2 \rangle = 0 \quad (n>0) \;, $$
   $$ a_0 | l_1, l_2 \rangle = l_1 | l_1, l_2 \rangle
      \; , \quad
      b_0 | l_1, l_2 \rangle = l_2 | l_1, l_2 \rangle \; , $$
   where $| l_1, l_2 \rangle$ is the vacuum vector.
   The grading operator $\overline{d}$ is defined by
   $$ \overline{d}.f |l_1,l_2 \rangle
        =(-\sum^k_{l=1} m_l-\sum^{k'}_{l=1} n_l+\frac{l_1^2-l_2^2+l_2}{2})
         f |l_1,l_2 \rangle $$
   for $ f=a_{-m_1} \cdots a_{-m_k}b_{-n_1} \cdots b_{-n_{k'}} $.
   We set 
   $$ \widetilde{\cal F}:=
          \bigoplus_{l_1 \in \frac{1}{2}{\bf Z},l_2 \in {\bf Z}}
          {\cal F}_{l_1,l_2} . $$
   We define the operators $e^{P_a}$ and $e^{P_b}$
   on ${\widetilde{\cal F}}$ as follows.
   $$ e^{P_a}|l_1,l_2 \rangle= |l_1+1,l_2 \rangle
      \quad , \quad
      e^{P_b}|l_1,l_2 \rangle= |l_1,l_2+1 \rangle . $$
   Let $: \quad :$ be the normal ordering defined by
   $$ :a_m a_n: = a_m a_n \, (m \leq n) , \; a_n a_m \, (m>n) , $$
   $$ :b_m b_n: = b_m b_n \, (m \leq n) , \; b_n b_m \, (m>n) , $$
   $$ :e^{P_a} a_0:=:a_0 e^{P_a}:=e^{P_a} a_0 \, , $$
   $$ :e^{P_b} b_0:=:b_0 e^{P_b}:=e^{P_b} b_0 \, . $$
   Consider the following currents.
   \[ Y^{\pm}_a(z)=
         \exp ( \pm \sum^{\infty}_{k=1}
                 \frac{a_{-k}}{[-\frac{1}{2}k]} q^{\pm \frac{k}{4}} z^k)
         \exp ( \mp \sum^{\infty}_{k=1}
                 \frac{a_k   }{[-\frac{1}{2}k]} q^{\pm \frac{k}{4}} z^{-k})
         e^{\pm 2P_a} z^{\mp 2a_0} , \]
   \[ Y^{\pm}_b(z)=
             \exp ( \pm \sum^{\infty}_{k=1} \frac{b_{-k}}{k} z^k)
             \exp ( \mp \sum^{\infty}_{k=1} \frac{b_k}{k} z^{-k})
             e^{\pm P_b} z^{\pm b_0} . \]

   \bigskip
   \noindent
   {\bf Proposition 3.1.}
   {\it
    $\widetilde{\cal F}$ is a $U_q$-module of level $-\frac{1}{2}$
    under the action of $U_q$ defined by 
    $$ \gamma \longmapsto q^{-\frac{1}{2}} \; \; \; ,
       \; \;
       K \longmapsto q^{a_0} \; \; , $$
    $$ h_n \longmapsto a_n \; \; \; ,
       \; \;
       q^d \longmapsto q^{\overline{d}} \; \; , $$
   \[ X^+(z) \longmapsto
             \left(
             \frac{1}{q^{\frac{1}{2}}+ q^{-\frac{1}{2}}}
             :Y^+_b(z) _{q^{\frac{1}{2}}} \partial^2_z Y^+_b(z):
             \right. \]
   \[   \hspace{1.5in}
             \left.
            -:_{q^{\frac{1}{2}}} \partial_z Y^+_b(q^{\frac{1}{2}}z)
              _{q^{\frac{1}{2}}} \partial_z Y^+_b(q^{-\frac{1}{2}}z):
             \right)
             Y^+_a(z) \; , \]
   \[ X^-(z) \longmapsto
             \frac{1}
                  {q^{\frac{1}{2}}+ q^{-\frac{1}{2}}}
             :Y^-_b(q^{\frac{1}{2}}z)Y^-_b(q^{-\frac{1}{2}}z): Y^-_a(z) \]

   where
   $$ X^{\pm}(z):=\sum_{k \in Z} x^{\pm}_k z^{-k-1} . $$
   }

   \bigskip
   The proof will be given in Section 5.1.

   The weight of $| l_1,l_2 \rangle$ is
   $ \frac{1}{2}l_1 \alpha_1 - \frac{1}{2}(l_1^2-l_2^2+l_2)\delta$.
   The space $\widetilde{\cal F}$ includes
   four irreducible highest weight submodules
   as we will see below.
   Consider the operator defined by
   \[ Q^-:=\oint Y^-_b(z) \frac{dz}{2 \pi \sqrt{-1}} \; : \;
      {\cal F}_{l_1,l_2} \longrightarrow {\cal F}_{l_1,l_2-1}. \]
   We define ${\cal F}_i$ $(i=1,2,3,4)$ by
   \[ {\cal F}_i 
      := \bigoplus_{l \in 2{\bf Z}}
         \Ker(Q^-: {\cal F}_{l+r_i,l+t_i}
                  \rightarrow
                  {\cal F}_{l+r_i,l+t_i-1}) . \]
   where $(r_i,t_i)=(0,0)$, $(1,1)$, $(-1/2,0)$, $(1/2,1)$
   for $i=1$, $2$, $3$, $4$, respectively.
   Then we have the following theorem.

   \bigskip
   \noindent
   {\bf Theorem 3.2.}
   {\it
   Each ${\cal F}_i$ $(i=1,2,3,4)$ is
   an irreducible highest weight $U_q$-module
   isomorphic to $V(\lambda'_i)$,
   where
   $\lambda_1' = - \frac{1}{2} \Lambda_0 $,
   $\lambda_2' = - \frac{3}{2} \Lambda_0 + \Lambda_1 - \frac{1}{2}\delta $,
   $\lambda_3' = - \frac{1}{2} \Lambda_1             + \frac{1}{8}\delta $,
   $\lambda_4' = \Lambda_0 - \frac{3}{2} \Lambda_1   + \frac{1}{8}\delta $.
   The highest weight vectors are given by
   $ | \lambda_1' \rangle =        |   0  ,  0 \rangle $,
   $ | \lambda_2' \rangle = b_{-1} |   1  ,  1 \rangle $,
   $ | \lambda_3' \rangle =        | -1/2 ,  0 \rangle $,
   $ | \lambda_4' \rangle =        | -3/2 , -1 \rangle $. }

   \bigskip
   \noindent
   {\it Proof.}
   It can be checked immediately
   that each  $| \lambda'_i \rangle$
   is a weight vector of weight $\lambda'_i$,
   satisfies the highest weight condition
   and belongs to ${\cal F}_i$.
   Since $[U_q, Q^-]=0$ (which will be proved in Section 5.2),
   ${\cal F}_i$ is a $U_q$-module.
   Next, we calculate the character of ${\cal F}_i$.
   We may understand $Q^-$ as the zero mode $\eta_0$
   of the fermionic ghost system $(\xi, \eta)$ of dimension $(0,1)$.
   $$ \xi(z) =Y^+_b(z)=\sum_{n \in {\bf Z}} \xi_n z^{-n} \; , \quad
      \eta(z)=Y^-_b(z)=\sum_{n \in {\bf Z}}\eta_n z^{-n-1} . $$
   Since we have $\eta_0^2=0$ and $\xi_0 \eta_0 + \eta_0 \xi_0 = 1$,
   we obtain the following exact sequence for any $l_1$, $l_2$
   $$ 0 \longrightarrow \Ker_{{\cal F}_{l_1,l_2}} Q^-
        \longrightarrow {\cal F}_{l_1,l_2}
        \stackrel{Q^-}{\longrightarrow} {\cal F}_{l_1,l_2-1}
        \stackrel{Q^-}{\longrightarrow} {\cal F}_{l_1,l_2-2}
        \stackrel{Q^-}{\longrightarrow} \cdots . $$
   Using this exact sequence, we can compute the character of ${\cal F}_i$.
   Since $a_n$ and $b_n$ have weight $n\delta$, we can easily see
   $$ \ch {\cal F}_{l_1,l_2}
      =\frac{z^{-\frac{1}{2} l_1} p^{-\frac{1}{2}(l_1^2-l_2^2+l_2)}}
            {(p;p)_{\infty}^2} , $$
   and hence
   \[ \ch({\cal F}_1 \oplus {\cal F}_2)
      =\sum_{m \in {\bf Z}} \ch (\Ker_{{\cal F}_{m,m}} Q^-) \]
   \[ \hspace{.91in}
      =(p;p)_{\infty}^{-2}
       \sum_{m \in {\bf Z}} z^{-\frac{1}{2}m}
             \sum_{k \leq m}
              (-1)^{m-k} p^{-\frac{1}{2}(m^2-k^2+k)} \]
   \[ \hspace{.91in}
       \stackrel{(*)}{=}
       \frac{1}{(p^{\frac{1}{2}} z^{ \frac{1}{2}};p)_{\infty}
                (p^{\frac{1}{2}} z^{-\frac{1}{2}};p)_{\infty}} . \]
   The last equality $(*)$ will be proved in Section 5.3.
   Similarly, we have
   $$ \ch({\cal F}_3 \oplus {\cal F}_4)
         =\frac{p^{-\frac{1}{8}} z^{\frac{1}{4}}}
               {(p z^{ \frac{1}{2}};p)_{\infty}
                (  z^{-\frac{1}{2}};p)_{\infty}} . $$
   Hence, by the character formulas in Section 2,
   $$ \ch{\cal F}_i = \ch L(\lambda_i') . $$
   As in Lusztig(\cite{l}),
   ${\cal F}_i$ becomes a certain $\widehat{sl}_2$-module
   in the classical limit $q \rightarrow 1$.
   Since the dimension of each weight space is invariant in the limit,
   the $\widehat{sl}_2$-module in the classical limit is irreducible.
   Therefore ${\cal F}_i$ is irreducible for a generic $q$.
   \hspace{\fill} $\Box$

   \bigskip
   We obtain $V(\lambda_i)$ from $V(\lambda'_i)$
   by shifting the grading operator $d$ to $d_i$,
   where $d_i=d$, $d+\frac{1}{2}$, $d-\frac{1}{8}$, $d-\frac{1}{8}$,
   for $i=1$, $2$, $3$, $4$, respectively.
   This representation becomes
   the representation constructed
   by Feingold and Frenkel in \cite{ff}
   in the classical limit $q \rightarrow 1$.
   The relation between $a_n$, $b_n$ and $\phi_1(z)$, $\phi_2(z)$
   in the introduction is as follows.
   $$ Y^{\pm}_a(z)
      \longrightarrow
      :e^{\pm 2 \phi_1(z)}: \; , \quad
     Y^{\pm}_b(z)
      \longrightarrow
      :e^{\pm   \phi_2(z)}: \; . $$


 \section{The $q$-vertex operators}

   \subsection{Definition of the $q$-vertex operators}

    We recall
    the definition and some properties
    of the $q$-vertex operators
    in our case(\cite{fr},\cite{ay}).
    We consider the 2-dimensional $U'_q$ module
    $ V = {\bf K}v_+ \oplus {\bf K}v_- $.
    The $U'_q$-module structure on $V$ is given by
    $$ e_1.v_+ = f_0.v_+ = 0,   \;
       e_1.v_- = f_0.v_- = v_+, \; $$
    $$ e_0.v_- = f_1.v_- = 0,   \;
       e_0.v_+ = f_1.v_+ = v_-, \; $$
    $$ t_0.v_{\pm}= q^{\mp 1} v_{\pm}, \;
       t_1.v_{\pm}= q^{\pm 1} v_{\pm} . $$
    The affinization of $V$
    is the following $U_q$-module $V_z$.
    $$ V_z := V \otimes {\bf K}[z,z^{-1}] . $$
    We define
    the $U_q$-module structure on $V_z$
    as follows.
    $$ e_i. (v \otimes z^m) = e_i. v \otimes z^{m+\delta_{i0}}
       \quad , \quad
       f_i. (v \otimes z^m) = f_i. v \otimes z^{m-\delta_{i0}}  $$
    $$ t_i. (v \otimes z^m) = t_i. v \otimes z^m
       \quad , \quad
       q^d. (v \otimes z^m) = q^m v \otimes z^m.     $$
    \\
    {\bf Definition 4.1.(\cite{fr})}
    {\it
    The $q$-vertex operator is
    a $U_q$-homomorphism of one of the following types.

    Type  I :
    $$ \tilde{\Phi}^{\mu V}_{\lambda}(z):
          V(\lambda) \longrightarrow V(\mu) \widehat{\otimes} V_z     $$

    Type II :
    $$ ^V\tilde{\Phi}^{\mu}_{\lambda}(z):
          V(\lambda) \longrightarrow V_z \widehat{\otimes} V(\mu)     $$}
    \hspace{\fill} $\Box$ 
    \\
    The existence conditions of the $q$-vertex operators
    are known in \cite{ay}.
    In our case, the $q$-vertex operators 
    exist only for
    $(\lambda, \mu)=(\lambda_1, \lambda_2)$,
    $(\lambda_2, \lambda_1)$,
    $(\lambda_3, \lambda_4)$,
    $(\lambda_4, \lambda_3)$.
    Furthermore, each of them is unique up to a scalar.
    Here, we take the following normalization.
    $$ \tilde{\Phi}^{\mu V}_{\lambda}(z)
       | \lambda \rangle
          =|\mu \rangle \otimes v_+
             + ( {\rm \; the \; terms \; of
                      \; positive \; powers \; in} \; z \; ) $$
    $$ \qquad \qquad \qquad \qquad
       {\rm for} \;
       (\lambda,\mu)=(\lambda_2,\lambda_1),(\lambda_3,\lambda_4) , $$
    $$ \tilde{\Phi}^{\mu V}_{\lambda}(z)
       | \lambda \rangle
          =|\mu \rangle \otimes v_-
             + ( {\rm \; the \; terms \; of
                      \; positive \; powers \; in} \; z \; ) $$
    $$ \qquad \qquad \qquad \qquad {\rm for} \;
       (\lambda,\mu)=(\lambda_1,\lambda_2),(\lambda_4,\lambda_3) . $$
    For the type II, we take a similar normalization.

   \subsection{Bosonizations}

    In this subsection,
    we construct an explicit form of the $q$-vertex operators
    on ${\cal F}_i$'s.
    We define
    the components $\tilde{\Phi}^{\mu V}_{\lambda}{}_{\pm} (z)$
    of the $q$-vertex operators as follows.
    $$ \tilde{\Phi}^{\mu V}_{\lambda}(z) |u \rangle
           =  \tilde{\Phi}^{\mu V}_{\lambda}{}_+ (z)
              |u \rangle \otimes v_+
             +\tilde{\Phi}^{\mu V}_{\lambda}{}_- (z)
              |u \rangle \otimes v_-
         \qquad {\rm for} \quad |u \rangle \in V(\lambda) , $$
    For the type II, the components are defined similarly.

    Using the same technique as in \cite{jmmn} (c.f. \cite{k},\cite{jkk}),
    the components of the $q$-vertex operators are given by the following
    theorem.

    \bigskip
    \noindent
    {\bf Theorem 4.2.}

    \[ 1) \quad
          \widetilde{\Phi}^{\mu V}_{\lambda}{}_-(z)
         ={\cal J}^+(q^{\frac{3}{2}}z)
          \, _{q^{\frac{1}{2}}}\partial_z Y^+_b(q^{\frac{3}{2}}z) \, r , \]
    \[ \qquad
       \tilde{\Phi}^{\mu V}_{\lambda}{}_+(z)
       =[ \; \tilde{\Phi}^{\mu V}_{\lambda}{}_-(z)
          \; , f_1 \; ]_q \; , \]
    \[ \qquad \quad
       r=1,
         -q^{\frac{3}{2}}z,
         -q^{\frac{3}{4}} z^{\frac{1}{2}},
         -q^{\frac{3}{4}} z^{\frac{1}{2}} \]
    \[ \qquad \quad
       for \;
       (\lambda,\mu)=(\lambda_1,\lambda_2),
                     (\lambda_2,\lambda_1),
                     (\lambda_3,\lambda_4),
                     (\lambda_4,\lambda_3), \]               
   \[ 2) \quad ^V \widetilde{\Phi}^{\mu}_{\lambda}{}_+(z)
        ={\cal J}^-(q^{-\frac{1}{2}}z)
         Y^-_b(q^{-\frac{1}{2}}z) \, r , \]
    \[ \qquad
       ^V \tilde{\Phi}^{\mu}_{\lambda}{}_{-}(z)
       =[ \; ^V \tilde{\Phi}^{\mu}_{\lambda}{}_{+}(z)
          \; , e_1 \; ]_q \; , \]
    \[ \qquad \quad
       r=-q^{-1},
          q^{-\frac{1}{2}}z,
          q^{-\frac{1}{4}} z^{\frac{1}{2}},
          q^{-\frac{5}{4}} z^{\frac{1}{2}} \]
    \[ \qquad \quad
       for \;
       (\lambda,\mu)=(\lambda_1,\lambda_2),
                     (\lambda_2,\lambda_1),
                     (\lambda_3,\lambda_4),
                     (\lambda_4,\lambda_3), \]
   {\it where}
   $$ [X,Y]_q=XY-qYX , $$
   $$ {\cal J}^{\pm}(z)
      =  \exp (\pm \sum^{\infty}_{k=1}
                  \widetilde{a}_{-k} q^{ \mp \frac{1}{4}k} z^k)
         \exp (\pm \sum^{\infty}_{k=1}
                  \widetilde{a}_k    q^{ \mp \frac{1}{4}k} z^{-k})
         e^{\pm P_a} z^{\mp a_0} , $$
   $$ \widetilde{a}_k
         =\frac{q^{\frac{1}{2}k}+q^{-\frac{1}{2}k}}{[2k]}a_k . $$
   \\
   {\it Proof.}
   The proof of the intertwining relations will be given
   in Section 5.4.
   It can be checked directly
   that the operators defined in $1)$ and $2)$
   anticommute with $Q^-$.
   Therefore, they give operators on the ${\cal F}_i$'s. 
   \hspace{\fill} $\Box$

   We give an example of
   the two-point functions of the $q$-vertex operators
   by using the above bosonization formulas.

   \[ \langle \lambda_1 |
      \widetilde{\Phi}^V (z_2)\widetilde{\Phi}^V (z_1)
      | \lambda_1 \rangle
      = \frac{(q^3 z;q^4)_{\infty}(q^6 z;q^4)_{\infty}}
             {(q   z;q^4)_{\infty}(q^4 z;q^4)_{\infty}}
        (v_+ \otimes v_- -q v_- \otimes v_+) , \]

   where $z=z_1/z_2$.

%
%
%
%


   \section{Proofs}

   In this section, we give the proofs postponed in the previous sections.
   We list the following formulas on 
   $Y^{\pm}_a(z)$, $Y^{\pm}_b(z)$ and ${\cal J}^{\pm}(z)$
   defined in the previous sections.

   \[ (1) \quad
       Y^{\pm}_a(z) Y^{\pm}_a(w)
      =\frac{:Y^{\pm}_a(z) Y^{\pm}_a(w):}
            {(z-q^{\pm 2}w)(z-qw)(z-w)(z-q^{-1}w)} , \]

   \[ (2) \quad
       Y^{\pm}_a(z) Y^{\mp}_a(w) \]
   $$ =:Y^{\pm}_a(z) Y^{\mp}_a(w):
        (z-q^{ \frac{3}{2}}w)
        (z-q^{ \frac{1}{2}}w)
        (z-q^{-\frac{1}{2}}w)
        (z-q^{-\frac{3}{2}}w) , $$

   \[ (3) \quad
       Y^{\pm}_b(z) Y^{\pm}_b(w)
      =:Y^{\pm}_b(z) Y^{\pm}_b(w):(z-w) , \]

   \[ (4) \quad
       Y^{\pm}_b(z) Y^{\mp}_b(w)
      =\frac{:Y^{\pm}_b(z) Y^{\mp}_b(w):}{z-w} , \]

   \[ (5) \quad
       {\cal J}^{\pm}(z) Y^{\pm}_a(w)
      =:{\cal J}^{\pm}(z) Y^{\pm}_a(w):
       \frac{1}{(z-q^{\frac{1}{2}}w)(z-q^{-\frac{1}{2}}w)} , \]

   \[ (6) \quad
       Y^{\pm}_a(w) {\cal J}^{\pm}(z)
      =:{\cal J}^{\pm}(z) Y^{\pm}_a(w):
       \frac{1}{(w-q^{\frac{1}{2}}z)(w-q^{-\frac{1}{2}}z)} , \]

   \[ (7) \quad
       {\cal J}^{\pm}(z) Y^{\mp}_a(w)
      =:{\cal J}^{\pm}(z) Y^{\mp}_a(w):
       (z-w)(z-q^{\mp 1}w) , \]

   \[ (8) \quad
       Y^{\mp}_a(w) {\cal J}^{\pm}(z)
      =:{\cal J}^{\pm}(z) Y^{\mp}_a(w):
       (w-z)(w-q^{\mp 1}z) , \]
   where $\frac{1}{z-w}$
   means the formal power series
   $\frac{1}{z} \sum_{k \geq 0} (\frac{w}{z})^k$.
   The above equations should be regarded
   as formal power series in $z$ and $w$.
   In the proofs, we will use the formal delta function
   $\delta(z)=\sum_{k \in {\bf Z}} z^k$.

   \subsection{Proof of Proposition 3.1.}

   Before starting the proof,
   we recall the defining relations
   of $U_q$ in terms of the generating functions
   $ X^{\pm}(z)=\sum_{k \in {\bf Z}} x^{\pm}_k z^{-k-1} $
   in the level $-1/2$ case.

     \[ R1) \hspace{.3in}
            K K^{-1} = K^{-1} K = 1 , \quad
            q^d q^{-d} = q^{-d} q^d = 1 , \quad
            K q^d = q^d K ,\]

     \[ R2) \hspace{.3in}
            q^d h_k q^{-d} = q^k h_k , \quad
            q^d X^{\pm}(z) q^{-d} = q^{-1} X^{\pm}(q^{-1}z) , \]

     \[ R3) \hspace{.3in}
            [ h_k , K ]=0 , \quad
            [ h_k , h_l ]
              =\delta_{k+l,0} \frac{[2k][-\frac{1}{2}k]}{k} , \]

     \[ R4) \hspace{.4in}
            K X^{\pm}(z) K^{-1}
                =q^{\pm 2} X^{\pm}(z) , \quad 
            [ a_k , X^{\pm}(z) ]
                =\pm \frac{[2 k]}{k}
                     q^{\pm \frac{|k|}{4}} X^{\pm}(z) , \]

     \[ R5) \hspace{.5in}
            (z-q^{-2}w)X^-(z)X^-(w)
           +(w-q^{-2}z)X^-(w)X^-(z) =0 \; , \]

     \[ R6) \hspace{.5in}
            (z-q^2 w)X^+(z)X^+(w)
           +(w-q^2 z)X^+(w)X^+(z) =0 \; , \]

     \[ R7) \hspace{.4in} [X^+(z),X^-(w)] \]
     \[    \hspace{.75in}
             =\frac{1}{q-q^{-1}}
                 \frac{1}{zw}
                 \left(
                 \delta(\frac{zq^{\frac{1}{2}}}{w})
                 \psi(zq^{\frac{1}{4}})
                -\delta(\frac{z}{wq^{\frac{1}{2}}})
                 \varphi(w^{-1}q^{-\frac{1}{4}})
                 \right) \]
           \hspace{.7in} {\it where}
     \[    \hspace{.5in}
               \begin{array}{l}
                      \psi (z)
                     =K \, \exp \left( (q-q^{-1})
                             \sum^{\infty}_{k=1} h_k z^{-k}
                                                       \right) \; ,
                 \\
                      \varphi (z)
                     =K^{-1} \exp \left( -(q-q^{-1})
                                  \sum^{\infty}_{k=1} h_{-k} z^{-k}
                                                       \right) \; .
               \end{array} \]

   The relations $R1)$, $R2)$ and $R3)$ are trivial.
   The relation $R4)$ follows from the following formulas.

  \[ q^{a_0} Y^{\pm}_a(z) q^{-a_0}=q^{\pm2} Y^{\pm}_a(z) , \]

  \[ [ a_k , Y^{\pm}_a(z) ]
     =\pm \frac{[2k]}{k} q^{\pm \frac{|k|}{4}} z^k Y^{\pm}_a(z)
     \qquad {\rm for} \ k \neq 0 . \]
   For the proof of relations $R5)$, $R6)$ and $R7)$,
   we introduce $M^+_1(z)$, $M^+_2(z)$, $M^+_3(z)$ and $M^{\pm}(z)$
   as follows.   

   \[ M^+_1(z) = : Y^+_a(z) Y^+_b(qz) Y^+_b(z) : , \]

   \[ M^+_2(z) = : Y^+_a(z) Y^+_b( z) Y^+_b(q^{-1}z) : , \]

   \[ M^+_3(z) = : Y^+_a(z) Y^+_b(qz) Y^+_b(q^{-1}z) : , \]

   \[ M^+(z) = q^{ \frac{1}{2}} M^+_1(z)
               +q^{-\frac{1}{2}} M^+_2(z)
               -(q^{\frac{1}{2}}+q^{-\frac{1}{2}})M^+_3(z) , \]

   \[ M^-(z) =:Y^+_a(z)Y^-_b(q^{\frac{1}{2}}z)Y^-_b(q^{-\frac{1}{2}}z): . \]
   The action of $X^{\pm}(z)$ is expressed by these currents as follows.
   \[ X^+(z) \longmapsto 
            -\frac{M^+(z)}{(q-q^{-1})(q^{\frac{1}{2}}-q^{-\frac{1}{2}})z^2}
      \; , \quad
      X^-(z) \longmapsto
             \frac{M^-(z)}{q^{\frac{1}{2}}+q^{-\frac{1}{2}}} . \]
   Using $(1)-(4)$, we see
   \[ (z-q^{-2}w)M^-(z)M^-(w)=:M^-(z)M^-(w):(z-w) . \]
   The relation $R5)$ follows from this.

   Next, we check the relation $R6)$.
   \[ \quad
      (z-q^2 w) M^+(z) M^+(w) \]
   \[ = :M^+_1(z) M^+_1(w): q^3 (z-w) \]
   \[ \quad
       +:M^+_1(z) M^+_2(w):
        q^2 \frac{(z-q^{-1} w)(z-q^{-2} w)}{z-qw} \]
   \[ \quad
       +:M^+_1(z) M^+_3(w):
        q^{\frac{5}{2}} (q^{\frac{1}{2}}+q^{-\frac{1}{2}})
        (z-q^{-2}w) \]
   \[ \quad
       +:M^+_2(z) M^+_1(w):
       q^{-2} \frac{(z-q w)(z-q^2 w)}{z-q^{-1}w} \]
   \[ \quad
       +:M^+_2(z) M^+_2(w):
       q^{-3} (z-w) \]
   \[ \quad
       +:M^+_2(z) M^+_3(w):
       q^{-\frac{5}{2}} (q^{\frac{1}{2}}+q^{-\frac{1}{2}})
       (z-q^2 w) \]
   \[ \quad
      +:M^+_3(z) M^+_1(w):
       q^{\frac{1}{2}} (q^{\frac{1}{2}}+q^{-\frac{1}{2}})
       (z-q^2 w) \]
   \[ \quad
      +:M^+_3(z) M^+_2(w):
       q^{-\frac{1}{2}} (q^{\frac{1}{2}}+q^{-\frac{1}{2}})
       (z-q^{-2} w) \]
   \[ \quad
      +:M^+_3(z) M^+_3(w):
       (q^{\frac{1}{2}}+q^{-\frac{1}{2}})^2
       \frac{(z-w)(z-q^2 w)(z-q^{-2} w)}
            {(z-qw)(z-q^{-1}w)} . \]
   Symmetrizing it with respect to $z$ and $w$,
   we have
   \[ \quad
       (z-q^2 w) M^+(z) M^+(w) + (w-q^2 z) M^+(w) M^+(z) \]
   \[ = :M^+_1(z) M^+_2(w):
        q (z-q^{-1}w)(z-q^{-2}w)
        \frac{1}{w} \delta (\frac{qw}{z}) \]
   \[ \qquad
     + :M^+_2(z) M^+_1(w):
       q^{-1} (z-qw)(z-q^2 w)
       \frac{1}{w} \delta (\frac{w}{qz}) \]
   \[ \qquad
      + :M^+_3(z) M^+_3(w): \]
   \[ \qquad \quad \times
             \frac{(q^{\frac{1}{2}}+q^{-\frac{1}{2}})^2
             (z-w)(z-q^2 w)(z-q^{-2} w)
             \left(
             \delta (\frac{qw}{z}) - \delta (\frac{qz}{w})
             \right)}
            {zw(q-q^{-1})} \]
   \[ = \Big(
        :M^+_1(qw) M^+_2(w):-:M^+_3(qw) M^+_3(w):
        \Big) \]
   \[ \hspace{2in} \times
        q     (q^{-1}-q)(q-q^{-2})  w \delta(\frac{qw}{z}) \]
   \[ \qquad
       +\Big(
        :M^+_2(q^{-1}w) M^+_1(w):-:M^+_3(q^{-1}w) M^+_3(w):
        \Big) \]
   \[ \hspace{2in} \times
        q^{-1}(q^{-1}-q)(q^{-1}-q^2)w \delta(\frac{w}{qz}) \]
   \[ =0 . \]
   The relation $R6)$ follows from this.

  Finally, we check the relation $R7)$.
  \[ \quad
     [M^+_1(z),M^-(w)] \]
  \[ = :M^+_1(z)M^-(w): q^{-2} \frac{z-q^{\frac{3}{2}}w}{z-q^{-\frac{1}{2}}w}
       \;
      -:M^+_1(z)M^-(w): \frac{w-q^{-\frac{3}{2}}z}{w-q^{\frac{1}{2}}z} \]
  \[ =:M^+_1(z)M^-(w):
      (z-q^{-\frac{3}{2}}w)\delta (\frac{z q^{\frac{1}{2}}}{w}) \]
  \[ =(q^{-2}-1)q^{a_0}
      \exp \Big((q-q^{-1}) \sum_{n=1}^{\infty} a_n(q^{\frac{1}{4}}z)^{-n}\Big)
      \delta (\frac{z q^{\frac{1}{2}}}{w}) . \]
  Similarly,
  \[ \quad
     [M^+_2(z),M^-(w)] \]
  \[     =(q^2-1)q^{-a_0}
      \exp \Big(-(q-q^{-1})
           \sum_{n=1}^{\infty} a_{-n} (q^{-\frac{1}{4}}w^{-1})^{-n} \Big)
      \delta (\frac{z}{w q^{\frac{1}{2}}}) , \]
  \[ \quad
     [M^+_3(z),M^-(w)]=0 . \]
  The relations $R7)$ follows from the above formulas.
  The proof is completed.

   \subsection{Proof of $[U_q,Q^-]=0$}

   It is trivial
   that $Q^-$ commutes with $\gamma^{\frac{1}{2}}$, $K$, $h_k$ and $x^-_k$.
   We will show $[X^+(z),Q^-]=0$.
   We need to calculate the commutator $[ M^+_i(z) , Y^-_b(w) ]$. 
   \[ [ M^+_1(z) , Y^-_b(w) ]
     = :M^+_1(z)Y^-_b(w):
       \frac{1}{zw(q-1)}
       \left(
       \delta(\frac{w}{z})-\delta(\frac{qz}{w})
       \right) . \]
   \[ [ M^+_2(z) , Y^-_b(w) ]
     = :M^+_2(z)Y^-_b(w):
       \frac{1}{zw(q^{-1}-1)}
       \left(
       \delta(\frac{w}{z})-\delta(\frac{z}{qw})
       \right) . \]
   \[ [ M^+_3(z) , Y^-_b(w) ]
     = :M^+_3(z)Y^-_b(w):
       \frac{1}{zw(q-q^{-1})}
       \left(
       \delta(\frac{qw}{z})-\delta(\frac{qz}{w})
       \right) . \]
   Using these formulas, we have
   \[ [ M^+(z) , Y^-_b(w) ]
     =Y^+_a(z)
      \left(
       Y^+_b(qz)      \Big(\delta(\frac{w}{z})-\delta(\frac{qw}{z}) \Big)
      \right. \]
   $$ \left.
      +Y^+_b(z)       \Big(\delta(\frac{qw}{z})-\delta(\frac{w}{qz}) \Big)
      +Y^+_b(q^{-1}z) \Big(\delta(\frac{w}{qz})-\delta(\frac{w}{z}) \Big)
      \right) $$
   Picking up the coefficient of $w^{-1}$, we have

   $$ [ X^+(z), Q^-]=0 . $$

   \subsection{Proof of the identity in the proof of theorem 3.2.}

   In this subsection,
   we prove the following identity $(*)$
   used in the calculation of the characters.
   $$ (*) \quad
       (p;p)_{\infty}^{-2}
       \sum_{m \in {\bf Z}} z^{-\frac{1}{2}m}
       \sum_{k \geq -m}
       (-1)^{k+m} p^{\frac{1}{2}(k^2-m^2+k)}
      =\frac{1}{(p^{\frac{1}{2}} z^{ \frac{1}{2}};p)_{\infty}
                (p^{\frac{1}{2}} z^{-\frac{1}{2}};p)_{\infty}} . $$
   In order to prove the above formula,
   it is sufficient to show
   $$ \sum_{n \in {\bf Z}} (-1)^n p^{\frac{1}{2} n^2} z^{\frac{1}{2}n}
      \sum_{m \in {\bf Z}}  z^{-\frac{1}{2}m}
      \sum_{k \geq -m} (-1)^{k+m} p^{\frac{1}{2}(k^2+k-m^2)}
     =(p;p)^3_{\infty} . $$
   The left hand side is rewritten as
   \[ \sum_{l \in {\bf Z}} (-1)^l p^{\frac{1}{2}l^2} z^{\frac{1}{2}l} 
      \sum_{m \in {\bf Z}} p^{ml}
      \sum_{k \geq -m} (-1)^k p^{\frac{1}{2}(k^2+k)} . \]
   Now, we show
   $ S_l=\sum_{m \in {\bf Z}} p^{ml}
         \sum_{k \geq -m} (-1)^k p^{\frac{1}{2}(k^2+k)}
        =(p;p)^3_{\infty} \delta_{l,0}$.
   \[ (1-p^l)S_l=S_l-\sum_{m \in {\bf Z}} p^{(m+1)l}
                     \sum_{k \geq -m} (-1)^k p^{\frac{1}{2}(k^2+k)} \]
   \[ \hspace{.72in}
                =S_l-\sum_{m \in {\bf Z}} p^{ml}
                     \sum_{k \geq -m+1} (-1)^k p^{\frac{1}{2}(k^2+k)} \]
   \[ \hspace{.72in}
                =\sum_{m \in {\bf Z}} p^{ml} (-1)^m p^{\frac{1}{2}(m^2-m)} \]
   \[ \hspace{.72in}
                =(p;p)_{\infty}(p^l;p)_{\infty}(p^{1-l};p)_{\infty} \]
   \[ \hspace{.72in}
                =0 . \]
   Hence, $S_l=0$ for $l \neq 0$.
   \[ S_0=\sum_{m \in {\bf Z}}
          \sum_{k \geq -m} (-1)^k p^{\frac{1}{2}(k^2+k)}
         =\sum_{m \in {\bf Z}}
          \sum_{k \geq |m|} (-1)^k p^{\frac{1}{2}(k^2+k)} \]
   \[ \hspace{.21in}
         =\sum_{n \geq 0} (-1)^n (2n+1) p^{\frac{1}{2}(n^2+n)}
         =(p;p)^3_{\infty} . \]
   The identity was proved.

   \subsection{Proof of Theorem 4.2.}

   We give the proof for type I.
   The intertwining conditions with the Chevalley generators
   are as follows.
     \[ V1) \qquad
           t_1 \widetilde{\Phi}^V_{\pm}(z) t_1^{-1}
         = q^{\mp 1} \widetilde{\Phi}^V_{\pm}(z) \]
     \[ V2) \qquad
           t_0 \widetilde{\Phi}^V_{\pm}(z) t_0^{-1}
         = q^{\pm 1} \widetilde{\Phi}^V_{\pm}(z) \]
     \[ V3) \qquad
           [ \widetilde{\Phi}^V_+(z) , e_0 ] = 0 \]
     \[ V4) \qquad
           [ \widetilde{\Phi}^V_-(z) , e_0 ] = z t_0 \widetilde{\Phi}^V_+(z) \]
     \[ V5) \qquad
           [ \widetilde{\Phi}^V_+(z) , e_1 ] = t_1 \widetilde{\Phi}^V_-(z) \]
     \[ V6) \qquad
           [ \widetilde{\Phi}^V_-(z) , e_1 ] = 0 \]
     \[ V7) \qquad
           [ \widetilde{\Phi}^V_+(z) , f_0 ]_q
           = z^{-1} \widetilde{\Phi}^V_-(z) \]
     \[ V8) \qquad
           [ \widetilde{\Phi}^V_-(z) , f_0 ]_{q^{-1}}=0 \]
     \[ V9) \qquad
           [ \widetilde{\Phi}^V_-(z) , f_1 ]_q = \widetilde{\Phi}^V_+(z) \]
     \[ V10) \hspace{.245in}
           [ \widetilde{\Phi}^V_+(z) , f_1 ]_{q^{-1}}=0 \]
   where $[X,Y]_{q^{\pm 1}}=XY-q^{\pm 1}YX$.
   \\
   These conditions are not independent.
   Actually, as we will see below, some of them follow from the others.
   We recall our construction of $\widetilde{\Phi}^V_{\pm}(z)$.
   \[ \widetilde{\Phi}^V_-(z)
         = {\cal J}^+(q^{\frac{3}{2}}z)
          _{q^{\frac{1}{2}}}\partial_z Y^+_b(q^{\frac{3}{2}}z) , \]
   and $\widetilde{\Phi}^V_+(z)$ is given by $V9)$.
   From the construction, $V1)$, $V2)$ and $V9)$
   are trivial.
   First, we shall reduce the conditions
   to $V3)$ and the following two equations.

   \[ (A) \qquad
         [\widetilde{\Phi}^V_-(z),X^+(w)]=0 . \]
   \[ (B) \qquad
         [\widetilde{\Phi}^V_-(z),X^-(w)]_q
         =[\widetilde{\Phi}^V_-(z),X^-(w)]_{q^{-1}}
          \left(\frac{w}{q^{\frac{1}{2}}z} \right) . \]
   It can be checked immediately that $V6)$ and $V8)$ follow from $(A)$
   and that $V4)$ follows from $(B)$.

    $V5)$ follows from $V1)$, $V6)$ and $V9)$.
    \[ \quad
        [ \widetilde{\Phi}^V_+(z) , e_1 ] \]
    \[ = [ [ \widetilde{\Phi}^V_-(z) , f_1 ]_q , e_1 ] \]
    \[ = [ [ \widetilde{\Phi}^V_-(z) , e_1 ] , f_1]_q
        -[  \widetilde{\Phi}^V_-(z) , [ e_1 , f_1 ] ]_q \]
    \[ =-[  \widetilde{\Phi}^V_-(z) , [ e_1 , f_1 ] ]_q \]
    \[ =-[  \widetilde{\Phi}^V_-(z) , \frac{t_1-t_1^{-1}}{q-q^{-1}} ]_q \]
    \[ = t_1 \widetilde{\Phi}^V_-(z) . \]

    $V7)$ follows from $V2)$, $V4)$ and $V8)$.
    \[ \quad
        [ \widetilde{\Phi}^V_+(z) , f_0 ]_q \]
    \[ = [ z^{-1} t_0^{-1} [ \widetilde{\Phi}^V_-(z) , e_0 ] , f_0 ]_q \]
    \[ = z^{-1} t_0^{-1}
         [ [ \widetilde{\Phi}^V_-(z) , e_0 ] , f_0 ]_{q^{-1}} \]
    \[ = z^{-1} t_0^{-1}
         ([ [ \widetilde{\Phi}^V_-(z) , f_0 ]_{q^{-1}} , e_0 ]_q
         +[  \widetilde{\Phi}^V_-(z) , [ e_0 , f_0 ] ]_{q^{-1}} ) \]
    \[ = z^{-1} t_0^{-1}
         [ \widetilde{\Phi}^V_-(z) , [ e_0 , f_0 ] ]_{q^{-1}} \]
    \[ = z^{-1} t_0^{-1}
         [ \widetilde{\Phi}^V_-(z) ,
           \frac{t_0-t_0^{-1}}{q-q^{-1}} ]_{q^{-1}} \]
    \[ = z^{-1} \widetilde{\Phi}^V_-(z) . \]

    $V10)$ follows from $V4)$, $V9)$ and $V3)$.
    \[ \quad
        [ \widetilde{\Phi}^V_+(z) , f_1 ]_{q^{-1}} \]
    \[ =[ z^{-1} t_0^{-1}[ \widetilde{\Phi}^V_-(z) , e_0 ] , f_1 ]_{q^{-1}} \]
    \[ = z^{-1} t_0^{-1}[ [ \widetilde{\Phi}^V_-(z) , e_0 ] , f_1 ]_q \]
    \[ = z^{-1} t_0^{-1}[ [ \widetilde{\Phi}^V_-(z) , f_1 ]_q , e_0 ] \]
    \[ = z^{-1} t_0^{-1}[ \widetilde{\Phi}^V_+(z) , e_0 ] \]
    \[ =0 . \]
    Therefore, it remains to prove $V3)$, $(A)$ and $(B)$.

    Frist, we prove $(B)$.
    Using $(4)-(8)$, we have
    \[ [ \widetilde{\Phi}^V_-(z) , X^-(w) ]_q \times q^{\frac{3}{2}}z
       = -:Y^-_b(q^2z){\cal J}^+(q^{\frac{3}{2}}z)M^-(w):
          \frac{w}{(w-q^{\frac{5}{2}}z)} \]
    \[ \hspace{2in}
         -:Y^-_b(q  z){\cal J}^+(q^{\frac{3}{2}}z)M^-(w):
         \frac{q^{-\frac{1}{2}}w}{(z-q^{-\frac{1}{2}}w)} , \]
    \[ [ \widetilde{\Phi}^V_-(z) , X^-(w) ]_{q^{-1}} \times q^{\frac{3}{2}}z
       = -:Y^-_b(q^2z){\cal J}^+(q^{\frac{3}{2}}z)M^-(w):
          \frac{q^{\frac{1}{2}}z}{(w-q^{\frac{5}{2}}z)} \]
    \[ \hspace{2in}
         -:Y^-_b(q  z){\cal J}^+(q^{\frac{3}{2}}z)X^-(w):
         \frac{z}{(z-q^{-\frac{1}{2}}w)} . \]
    $(B)$ follows from these two formulas.

    Next, we prove $(A)$.
    \[ \quad
       [\widetilde{\Phi}^V_-(z),M^+_1(w)] \]
    \[ =-:{\cal J}^+(q^{\frac{3}{2}}z) Y^+_b(qz) M^+_1(w):
         \frac{q(z-w)}{(q^2-q)zw}
         \delta(\frac{q^2z}{w}) \]
    \[ = :{\cal J}^+(q^{\frac{3}{2}}z) Y^+_b(qz) M^+_1(q^2 z):
        (q+1) w^{-1} \delta(\frac{q^2 z}{w}) . \]
    Similarly,
    \[ \quad
       [\widetilde{\Phi}^V_-(z),M^+_2(w)] \]
    \[ = :{\cal J}^+(q^{\frac{3}{2}}z) Y^+_b(q^2z) M^+_2(qz):
         (q+1) w^{-1} \delta(\frac{qz}{w}) , \]
    \[ \quad
       [\widetilde{\Phi}^V_-(z),M^+_3(w)] \]
    \[ = :{\cal J}^+(q^{\frac{3}{2}}z) Y^+_b(q^2z) M^+_3(q^2z):
         qw^{-1} \delta(\frac{q^2 z}{w}) \]
    \[ \qquad \qquad
        +:{\cal J}^+(q^{\frac{3}{2}}z) Y^+_b(q  z) M^+_3(q  z):
          w^{-1} \delta(\frac{q z}{w}) . \]
    From the above three formulas, we have
    \[
       [\widetilde{\Phi}^V_-(z),M^+(w)]=0 . \]
    We proved $(A)$.

    Finally, we prove $V3)$.
    \[  M^-(w_1)M^-(w_2)\widetilde{\Phi}^V_-(z) \]
    \[ -(q+q^{-1})
        M^-(w_1)\widetilde{\Phi}^V_-(z)M^-(w_2)
       +\widetilde{\Phi}^V_-(z)M^-(w_1)M^-(w_2) \]
    \[ =:{\cal J}^+(q^{\frac{3}{2}}z) Y^-_b(q^2 z) M^-(w_1) M^-(w_2):
        \frac{q^{-4}(q^{\frac{1}{2}}+q^{-\frac{1}{2}})(w_1-w_2)}
             {(w_1-q^{\frac{5}{2}}z)(w_2-q^{\frac{5}{2}}z)} \]
    \[ \qquad
       -:{\cal J}^+(q^{\frac{3}{2}}z) Y^-_b(q^2 z) M^-(w_1) M^-(w_2):
        \frac{q^{-1}(q^{\frac{1}{2}}+q^{-\frac{1}{2}})(w_1-w_2)}
             {(z-q^{-\frac{1}{2}}w_1)(z-q^{-\frac{1}{2}}w_2)} . \]
    Since $:M^-(w_1) M^-(w_2):=:M^-(w_2) M^-(w_1):$, the above equation
    is antisymmetric with respect to $w_1$ and $w_2$.
    Hence, the coefficient of $w_1^{-1} w_2^{-1}$ is zero.
    We proved $V3)$.

    The proof for the type I $q$-vertex operators is completed.
    For type II, the proof is similar.

    {\footnotesize
    \bigskip
    \noindent
    {\it Acknowledgement.}
    The author would like to thank
    Prof. M. Kashiwara, Prof. T. Miwa and Prof. M. Jimbo 
    for their encouragement and valuable discussions. 
    He also thanks
    K.Iohara, T.Miyazaki, Y.Saito, T.Suzuki and Yu.Yamada
    for stimulating discussions.}



\end{document}